\documentclass[12pt]{amsart}
\usepackage{amssymb,latexsym,eufrak,amsmath,amscd}
\usepackage[all]{xy}
\theoremstyle{plain}
\newtheorem{theorem}{Theorem}[section]
\newtheorem{proposition}[theorem]{Proposition}
\newtheorem{corollary}[theorem]{Corollary}
\newtheorem{lemma}[theorem]{Lemma}
\newtheorem{variant}[theorem]{Variant}

\newtheorem*{theoremnn}{Theorem}
\newtheorem*{restrthm}{Restriction Theorem}
\newtheorem*{nadelvan}{Nadel Vanishing Theorem}

\theoremstyle{definition}
\newtheorem{definition}[theorem]{Definition}
\newtheorem{remark}[theorem]{Remark}

\newcommand{\lra}{\longrightarrow}
\newcommand{\noi}{\noindent}
\newcommand{\PP}{\mathbb{P}}

\newcommand{\NN}{\mathbb{N}}
\newcommand{\ZZ}{\mathbb{Z}}
\newcommand{\CC}{\mathbb{C}}
\newcommand{\QQ}{\mathbb{Q}}

\newcommand{\OO}{\mathcal{O}}
\newcommand{\JJ}{\mathcal{J}}
\newcommand{\cH}{\mathcal{H}}

\newcommand{\FF}{\mathcal{F}}

\newcommand{\fra}{\frak{a}}
\newcommand{\frb}{\frak{b}}
\newcommand{\eps}{\epsilon}

\newcommand{\HH}[3]{H^{{#1}} \big( {#2} , {#3} \big) }
\newcommand{\hh}[3]{h^{{#1}} \big( {#2} , {#3} \big) }

\newcommand{\for}{ \ \ \text{ for } \ }
\newcommand{\fall}{ \ \ \text{ for all } \ }

\newcommand{\exc}{\text{Exc}}
\newcommand{\MI}[1]{\mathcal{J} ( {#1} ) }
\newcommand{\MIP}[1]{\mathcal{J}_+ ( {#1} ) }

\newcommand{\pr}{\prime}

\newcommand{\lin}{\equiv}

\title {A Subadditivity Property of Multiplier Ideals}

\author{Jean-Pierre Demailly} 
\address{ Laboratoire de Math\'ematiques, URA 5582 du
CNRS \\ Universit\'e de Grenoble I \hfil\break\indent
38402 Saint-Martin d'H\`eres, France}
\thanks{Research of first author partially supported by
CNRS}
\email{demailly@ujf-grenoble.fr}

\author{Lawrence Ein}
\address{Department of Mathematics \\University of
Illinois at Chicago \hfil\break\indent 
851 South Morgan Street (M/C 249)\\ Chicago, IL 60607-7045, USA}
\email{ein@math.uic.edu}
\thanks{Research of the second author partially
supported by NSF Grant DMS  99-70295 }

\author{Robert Lazarsfeld} 
\address{Department of Mathematics
\\ University of Michigan \\ Ann Arbor, MI 48109, USA}
\email{rlaz@math.lsa.umich.edu}
\thanks{Research of  third author  partially supported
by the J. S. Guggenheim Foundation and NSF Grant DMS
97-13149 }

\dedicatory{Dedicated to William Fulton on the occasion
of his sixtieth birthday}

\begin{document}

\maketitle
\setcounter{section}{-1}


\section*{Introduction}

The purpose of this note is to establish a
``subadditivity" theorem for multiplier ideals. As an
application, we give a new proof of a theorem of 
Fujita concerning the volume of a big line bundle.

Let $X$ be a smooth  complex quasi-projective variety,
and let $D$ be an effective $\QQ$-divisor on $X$. One
can associate to $D$ its \textit{multiplier ideal} sheaf
\[ \MI{D} = \MI{X,D} \subseteq \OO_X, \]
 whose zeroes are supported on the locus at which the
pair $(X,D)$ fails to have log-terminal singularities.
 It is useful to think of $\MI{D}$ as reflecting in a
somewhat subtle way the singularities of $D$: the
``worse" the singularities, the smaller the ideal. 
These ideals and their variants have  come to play an
increasingly important role in higher dimensional
geometry, largely because of their strong vanishing
properties. Among the papers  in which they figure
prominently, we might mention for instance 
\cite{Nadel.89}, \cite{Demailly.93},  \cite{Siu.93},
\cite{Angehrn.Siu}, 
\cite{Ein.Laz.1}, \cite{Siu.97}, \cite{Kawamata},
\cite{Ein.Laz.2} and \cite{DemaillyKollar.99}.  See
\cite{Demailly.Bourbaki} for a survey. 

We establish  the following ``subadditivity" property
of these ideals:
\begin{theoremnn} Given any two effective $\QQ$-divisors
$D_1$ and $D_2$ on $X$, one has the relation
\[ \MI{D_1 + D_2} \ \subseteq \ \MI{D_1} \cdot \MI{D_2}.
\]
\end{theoremnn}
\noi The Theorem admits several variants. In the local
setting, one can associate a multiplier ideal
$\MI{\fra}$ to any ideal $\fra \subseteq \OO_X$, which
in effect measures the singularities of the divisor of a
general element of $\fra$. Then the statement becomes
\[ \MI{ \fra \cdot \frak{b} } \ \subseteq \ \MI{\fra}
\cdot
\MI{\frak{b}}.  \] On  the other hand, suppose that $X$
is a smooth projective variety, and $L$ is a big line
bundle on
$X$. Then one can define an ``asymptotic multiplier
ideal" $\MI{\Vert L \Vert} \subseteq \OO_X$, which
reflects the asymptotic behavior of the base-loci of
the linear series $|kL|$ for large $k$. In this setting
the Theorem shows that
\[ \MI{\Vert mL \Vert } \ \subseteq \ \MI{ \Vert L
\Vert}^{
\, m} . \] Finally, there is an analytic analogue
(which in fact implies the other statements): one can
attach a multiplier ideal to any plurisubharmonic
function on $X$, and then 
\[\MI{\phi + \psi} \subseteq \MI{\phi} \cdot
\MI{\psi} \] for any two such functions
$\phi$ and $\psi$.  The Theorem was suggested by a
somewhat weaker statement established  in
\cite{Demailly.99}. 

We apply the subadditivity relation to give a new proof
of a theorem of Fujita \cite{Fujita}. Consider 
a smooth projective variety $X$ of dimension $n$, and a
big line bundle $L$ on $X$.  The
\textit{volume} of $L$ is defined to be the positive
real number
\[  v(L) \ = \ \limsup_{k \to \infty} \, 
\frac{n!}{k^n} 
 \, \hh{0}{X}{\OO(kL)}. \] If $L$ is ample then $v(L) =
\int_X c_1(L)^n$, and in general (as we shall see) it
measures asymptotically the top self-intersection of
the ``moving part" of $|kL|$ (Proposition
\ref{moving.self.int}).  Fujita has established the
following
\begin{theoremnn} [Fujita, \cite{Fujita}] Given any
$\eps > 0$, there exists a birational modification \[
\mu : X^\pr = X^{\pr}_{\epsilon}  \lra X \]
 and a decomposition $\mu^* L \equiv E_{\epsilon} +
A_{\epsilon}$, where
$E = E_{\epsilon}$ is an effective $\QQ$-divisor and $A
= A_{\epsilon}$ an  ample
$\QQ$-divisor, such that $\big( A^n \big) > v(L) -
\epsilon$. \end{theoremnn}
\noi This would be clear if $L$ admitted a Zariski
decomposition, and so one thinks of the statement as a
numerical analogue of such a decomposition. Fujita's
proof of the Theorem is quite short, but rather tricky.
We give a new proof using multiplier ideals  which (to
the present authors at least) seems perhaps more
transparent. An outline of this approach to Fujita's
theorem appears also in
\cite{Demailly.99}.  We hope that these ideas may find
other applications in the future. 

The paper is divided into three sections. In the first,
we review (largely without proof)  the theory of
multiplier ideals from the algebro-geometric point of
view, and we discuss the connections between 
asymptotic algebraic constructions and their analytic
counterparts. The subadditivity theorem is established
in
\S 2, via an elementary argument using a ``diagonal"
trick as in
\cite{DemaillyKollar.99}. The application to Fujita's
theorem appears in \S 3, where as a corollary we deduce
a geometric description of the volume of a big line
bundle. 

We thank E. Mouroukos for valuable discussions. We are
especially delighted to have the opportunity to dedicate
this paper to William Fulton on the occasion of his
sixtieth birthday. His many contributions have done
much to shape contemporary algebraic geometry.  The
third author in particular --- having been first a
student and being now a colleague of Bill's --- has
learned a great deal  from Fulton  over the years. 

\section{Notation and Conventions}

\noi (0.1).  We work throughout with non-singular
algebraic varieties defined over the complex numbers
$\CC$.  

\noi (0.2).  We generally speaking do not distinguish
between line bundles and (linear equivalence classes
of) integral divisors. In particular, given a line
bundle
$L$, we write $\OO_X(L)$ for the corresponding
invertible sheaf on $X$, and we use additive notation
for the tensor product of line bundles. When $X$ is a
smooth variety, $K_X$ denotes as usual the canonical
divisor (class) on $X$.

\noi (0.3).  We write $\lin$ for linear equivalence of
$\QQ$-divisors: two such divisors $D_1, D_2$
are linear equivalent if and only if there is a
non-zero integer $m$ such that $mD_1 \lin mD_2$ in the
usual sense.

\section{Multiplier Ideals}

In this section we review  the construction and basic
properties of multiplier ideals from an
algebro-geometric perspective. For the most part we do
not give proofs; most can be found in
\cite{EsnaultViehweg} (Chapter 7), \cite{Ein.SC},
\cite{Ein.2} and
\cite{Kawamata}, and a detailed exposition will appear
in the forthcoming book \cite{PAG}.  The algebraic
theory closely parallels the analytic one, for which
the reader may consult \cite{Demailly.CIME}. We also
discuss in some detail the relationship between the
algebraically defined asymptotic multiplier ideals
$\MI{\Vert L \Vert}$ associated to a complete linear
series and their analytic counterparts. 

Let $X$ be a smooth complex quasi-projective variety,
and $D$ an effective $\QQ$-divisor on $X$. Recall that
 a \textit{log resolution} of $(X, D)$ is a proper
birational mapping  \[ \mu : X^\pr \lra X \] from a
smooth variety $X^\pr$ to $X$ having the property that
$\mu^* D + \exc(\mu)$ has simple normal crossing
support, $\exc(\mu)$   being the sum of the exceptional
divisors of $\mu$.
\begin{definition} \label{Def.Mult.Ideal} The
\textbf{multiplier ideal} of
$D$ is defined to be 
\begin{equation} \label{Acyclic.LB} \MI{D} \ = \
\MI{X,D} \ = \ \mu_*
  \OO_{X^\pr}
\big( K_{X^{\pr}/X} - [\mu^* D] \big).
\end{equation}
\end{definition}
\noi Here $K_{{X^\pr}/X}$ denotes the relative canonical
divisor $K_{{X^\pr}} - \mu^* K_X$, and as usual $[F]$
is the integer part or round-down of a $\QQ$-divisor
$F$. That $\MI{D}$ is indeed an ideal sheaf follows
from the observation that
$ \MI{D} \ \subseteq  \ \mu_*
\OO_{X^\pr}\big(K_{X^{\pr}/X}\big) \ = \ \OO_X$.  An
important point is that this definition is independent
of the choice of resolution. This can be verified
directly, but it also follows from the fact that
$\MI{D}$ has an analytic interpretation. 

Using the same notation as in \cite{Demailly.99}, we
take a plurisubharmonic function $\phi$ and denote by
$\MI{\phi}$ the sheaf of germs of holomorphic functions
$f$ on $X$ such that $\int |f|^2e^{-2\phi}dV$ converges
on a neighborhood of the given point. By a well-known
result of Nadel \cite{Nadel.89}, $\MI{\phi}$ is always
a coherent analytic sheaf,  whatever the singularities
of $\phi$ might be. In fact, this follows from 
H\"ormander's $L^2$ estimates (\cite{Kodaira.53}, \cite{Hormander.66},
\cite{Andreotti.Vesentini}) for the $\overline\partial$ operator, 
combined with some elementary arguments of local algebra 
(Artin-Rees lemma).  We need here a slightly more precise statement
which can be inferred directly from the proof given in \cite{Nadel.89}
(see also \cite{Demailly.93}).

\begin{proposition} Let $\phi$ be a plurisubharmonic
function on a complex manifold $X$, and let $U\subseteq
X$ be a relatively compact Stein  open subset $($with a
basis of Stein neighborhoods of~$\overline U)$. Then
the  restriction $\MI{\phi}_{|U}$ is generated as an
$\OO_U$-module by a  Hilbert basis $(f_k)_{k\in\NN}$ of
the Hilbert space $\cH^2(U,\phi, dV)$ of  holomorphic
functions $f$ on $U$ such that
\[
\int_U |f|^2e^{-2\phi}dV < +\infty
\] $($with respect to any K\"ahler volume form $dV$ on a
neighborhood of 
$\overline U\,)$. \qed
\end{proposition}

Returning to the case of an effective $\QQ$-divisor
$D = \sum a_i D_i$, let $g_i$ be a local defining
equation for $D_i$. Then, if $\phi$ denotes the
plurisubharmonic function $\phi =\sum a_i\log|g_i| $,
one has
\[
\MI{D} = \MI{ \phi},
\] and in particular $\MI{D}$ is intrinsically defined.
The stated equality is established in
\cite{Demailly.CIME}, (5.9): the essential point is that
the algebro-geometric multiplier ideals satisfy the
same transformation rule under birational
modifications as do their analytic counterparts, so
that one is reduced to the case where $D$ has normal
crossing support.  

 We mention two variants. First, suppose given an ideal
sheaf $\fra \subseteq \OO_X$. By a log resolution of
$\fra$ we understand a mapping $\mu : X^{\pr}
\lra X$ as above with the property that
$\fra \cdot \OO_{X^{\pr}} =
\OO_{X^\pr}(-E)$, where $E + \exc(\mu)$ has simple
normal crossing support. Given a rational number $c >
0$ we take such a resolution and then define   
\[ \MI{c \cdot \fra } \ = \ \mu_* \OO_{X^\pr}\big(
K_{X^{\pr}/X} - [c E] \big); \] again this is
independent of the choice of resolution.
\footnote{
 More generally, given ideals $\fra, \frb \subseteq
\OO_X$, and rational numbers $c , d > 0$, one can define
$\MI{ (c \cdot \fra ) \cdot (d \cdot \frb)} $ by taking
a common log resolution $\mu : X^{\pr} \lra X$ of
$\fra$ and $\frb$,  with $\mu^{-1} \fra =
\OO_{X^\pr}(-E_1)$ and
$\mu^{-1} \frb = \OO_{X^\pr}(-E_2)$, and setting 
\[ 
\MI{ (c \cdot \fra ) \cdot (d \cdot \frb)} = \mu_*
\OO_{X^{\pr}} \big( K_{X^{\pr}/X} - [ cE_1 + dE_2]
\big). \] } If $m \in \ZZ$ is a positive integer then
$\MI{m \cdot
\fra} = \MI{ \fra ^ m }$, and in this case these
multiplier ideals were defined and studied in a more
general setting by Lipman
\cite{Lipman} (who calls them ``adjoint ideals"). They
admit the following geometric interpretation. Working
locally, assume that $X$ is affine, view
$\fra$ as an ideal in its coordinate ring,  and take $k
> c$ general
$\CC$-linear combinations of a set of generators 
$g_1,  \ldots , g_p \in \fra$,  yielding divisors $A_1,
\ldots, A_k  \subset X$. If $D = \tfrac{c}{k} \big( A_1
+ \ldots + A_k \big)$, then
\begin{equation} \label{geom.interp} \MI{c \cdot \fra} \
= \ \MI{D}.
\end{equation} In the analytic setting, where $X$ is an
open subset of
$\CC^n $, one has $\MI{ c \cdot \fra } = \MI{ c
\cdot \phi}$, where $\phi =  \log ( |g_1| + \dots +
|g_p|) $.

The second variant involves linear series. Suppose that
$L$ is a line bundle on $X$, and that $V \subset
\HH{0}{X}{L}$ is a finite dimensional vector space of
sections of $L$, giving rise to a linear series
$|V|$ of divisors on $X$. We now require of our log
resolution $\mu : X^{\pr } \lra X$ that 
\[ \mu^* | V | \ = \ |W| + E, \]
 where $|W|$ is a free linear series on $X^\pr$, and $E
+ \exc(\mu)$ has simple normal  crossing support. In
other words, we ask that the fixed locus of $\mu^* |V|$
be a divisor $E$ with simple normal crossing  support
(which in addition meets
$\exc(\mu)$ nicely). Given such a log resolution, plus
a rational number $c > 0$  we define
\[ \MI{ \, c \cdot | V | \, } \ = \ \mu_*
\OO_{X^\pr} \big( K_{X^{\pr}/X} - [c E] \big), \] this
once again being independent of the choice of
$\mu$. If $ {\frak b} = {\frak  b} \big( |V| \big)
\subseteq
\OO_X$ is the base-ideal of $|V|$, then evidently
$\MI{c \cdot |V|} = \MI{ c \cdot {\frak b} }$, and in
particular the analogue of Equation (\ref{geom.interp})
holds for these ideals. 

We now outline the main properties of these ideals that
we shall require. The first is a local statement
comparing a multiplier ideal with its restriction to a
hyperplane. Specifically, consider an effective
$\QQ$-divisor $D$ on a quasi-projective complex manifold
$X$, and a smooth effective divisor $H \subset X$ which
does not appear in the support of $D$. Then one can
form two ideals on $H$. In the first place, the
restriction $D_{|H}$ is an effective $\QQ$-divisor on
$H$, and so one can form its multiplier ideal
$\MI{H, D_{|H}} \subseteq
\OO_H$. On the other hand, one can take the multiplier
ideal $\MI{X,D}$ of $D$ on $X$ and restrict it to $H$
to get an ideal 
\[ \MI{X,D} \cdot \OO_H \  \subseteq \ \OO_H . \] A
very basic fact --- due in the algebro-geometric
setting to Esnault-Viehweg \cite{EsnaultViehweg} --- is
that one can compare these sheaves:
\begin{restrthm}  In the setting just described, there
is an inclusion
\[ \MI{H, D_{|H}} \  \subseteq \ \MI{X,D} \cdot \OO_H.
\]
\end{restrthm}
\noi One may think of this as asserting that
``multiplier ideals can only get worse" upon
restricting a divisor to a hyperplane. For the proof,
see \cite{EsnaultViehweg}, (7.5), or \cite{Ein.SC},
(2.1).  The essential point is that the line bundle 
$ \OO_{X^\pr}
\big( K_{X^{\pr}/X} - [\mu^* D] \big)$ appearing in
Equation (\ref{Acyclic.LB}) has vanishing higher direct
images under $\mu$. The same result holds true in the
analytic case, namely
\[ \MI{H, \phi_{|H}} \  \subseteq \ \MI{X,\phi} \cdot
\OO_H \] for every plurisubharmonic function $\phi$ on
$X$ (if $\phi_{|H}$ happens to be identically equal to $-\infty$ on
some component of $H$, one agrees that $\MI{H, \phi_{|H}}$
is identically zero on that component). In that case, 
the proof is completely different; it is in fact a direct
qualitative consequence of  the (deep) Ohsawa-Takegoshi 
$L^2$ extension theorem \cite{OhsawaTakegoshi}, \cite{Ohsawa.88}.

As a immediate consequence, one obtains an analogous
statement for restrictions to submanifolds of higher
codimension:
\begin{corollary} \label{restr.to.submnfld} Let $Y
\subset X$ be a smooth subvariety which is not contained
in the support of $D$. Then
\[ \MI{Y, D_{|Y}} \ \subseteq \ \MI{X,D} \cdot \OO_Y ,
\] where $D_{|Y}$ denotes the restriction of $D$ to $Y$.
\qed
\end{corollary}
\noi Of course the analogous statement is still true in
the analytic case, as well as for the multiplier
ideals associated to linear series  or ideal sheaves. 

The most important global property of multiplier ideals
is the
\begin{nadelvan} Let $X$ be a smooth complex projective
variety, $D$ an effective $\QQ$-divisor and  $L$ a line
bundle on $X$. Assume that $L-D$ is big and nef. Then
\[ \HH{i}{X}{\OO_X(K_X + L) \otimes \MI{D}} = 0 \for i
> 0. \]
\end{nadelvan}
\noi This follows quickly from the Kawamata-Viehweg
vanishing theorem applied on a log resolution $\mu :
X^\pr \lra X$ of $(X,D)$. Similarly, if $V \subset
\HH{0}{X}{B}$ is a linear series on $X$, with $B$ a
line bundle such that $L - c \cdot B$ is big and nef, 
then
\[ \HH{i}{X}{\OO_X(K_X + L) \otimes \MI{c \cdot |V| }} =
0 \for i > 0. \] Under the same hypotheses, if $\fra
\subseteq \OO_X$ is an ideal sheaf such that $B \otimes
\fra$ is globally generated, then 
$ \HH{i}{X}{\OO_X(K_X + L) \otimes \MI{c \cdot \fra}} =
0$ when $i > 0$.

Nadel Vanishing  yields a simple criterion for a
multiplier ideal sheaf to be globally generated. The
essential point is the following elementary lemma of
Mumford, which forms the basis of the theory of
 Castelnuovo-Mumford regularity:\footnote{We beg the
reader's indulgence for the fact that we prefer to
state the Lemma using multiplicative notation for tensor
products of line bundles, rather than working
additively as we do elsewhere in the paper.}
\begin{lemma} [\cite{Mumford}, Lecture 14]
\label{CM.Reg} Let
$X$ be a projective variety, $B$ a very ample line
bundle on $X$, and $\FF$ any coherent sheaf on $X$
satisfying the vanishing
\[ \HH{i}{X}{\FF \otimes B^{\otimes (k - i)}} = 0 \for 
i > 0 \ \text{and} \  k \ge 0. \] Then $\FF$ is globally
generated. 
\end{lemma}
\noi Although the Lemma is quite standard, it seems not
to be as well known as one might expect in connection
with vanishing theorems (Remark
\ref{Remark.on.Reg}). Therefore we feel it is
worthwhile to write out the argument.
\begin{proof} Evaluation of sections determines a
surjective map
$e : H^0(B)\otimes_{\CC} \OO_X \lra B$ of vector bundles
on
$X$. The corresponding Koszul complex takes the form:
\begin{equation} \ldots\lra \Lambda^3 H^0(B) \otimes
B^{\otimes -2} 
\lra \Lambda^2 H^0(B)
\otimes B^{\otimes -1}
\lra H^0(B) \otimes \OO_X \lra B \lra 0. \tag{*}
\end{equation} Tensoring through   by $\FF$, and
applying the hypothesis with $k = 0$ as one chases
through the resulting complex, one sees first of all
that the multiplication map
\[ H^0(B) \otimes H^0(\FF)
\lra H^0(\FF
\otimes B)\]  is surjective. Next tensor (*)   by $\FF
\otimes B$ and apply the vanishing hypothesis with $k =
1$:  it follows that $H^0(B) \otimes H^0(\FF \otimes
B)$ maps onto $H^{0}(\FF \otimes B^{\otimes 2})$, and
hence that 
$H^0(B^{\otimes 2})
\otimes H^0(\FF)
\lra H^0(\FF \otimes B^{\otimes 2})$ is also onto.
Continuing,  one finds that 
\begin{equation} \HH{0}{X}{\FF} \otimes
\HH{0}{X}{B^{\otimes m}}
\lra
\HH{0}{X}{\FF \otimes B^{\otimes m}} \tag{**}
\end{equation}  is surjective for all $m \ge 0$. But
since
$B$ is very ample, $\FF \otimes B^{\otimes m}$ is
globally generated for
$m \gg 0$. It then follows from  the surjectivity of
(**) that
$\FF$ itself must already be generated by its global
sections. \footnote{A similar argument shows that the
case $k = 0$ of the vanishing hypothesis actually
implies the cases $k \ge 1$, but for present purposes
we don't need this.}
\end{proof}
\begin{corollary} \label{MI.Gl.Gen} In the setting of
the Nadel Vanishing Theorem, let $B$ be a very ample
line bundle on $X$. Then
\[ \OO_X(K_X + L + mB) \otimes \MI{D} \] is globally
generated for all $m \ge \dim X$.
\end{corollary}
\begin{proof} In fact, thanks to Nadel vanishing, the
hypothesis of Mumford's Lemma applies to $\FF =
\OO_X(K_X + L + mB)
\otimes \MI{D}$ as soon as $m \ge \dim X$. \end{proof}

\begin{remark} \label{Remark.on.Reg}  The Corollary was
used by Siu in the course of his spectacular proof   of
the deformation invariance of plurigenera
\cite{Siu.97}, where the statement was established by
analytic methods. Analogous applications of the Lemma in
the context of vanishing theorems have appeared
implicitly or explicitly in a number of papers over the
years, for instance
\cite{Wilson}, \cite{Kollar1},
\cite{EsnaultViehweg}, \cite{SAD} (to name a few). \qed
\end{remark}

We next turn to the construction of the asymptotic
multiplier ideal associated to a big linear series. In
the algebro-geometric setting, the theory is due to the
second author
\cite{Ein.1} and Kawamata \cite{Kawamata}. Suppose  that
$X$ is a smooth complex projective variety, and $L$ is
a big line bundle on $X$. Then $\HH{0}{X}{\OO_X(kL)}
\ne 0$ for $k \gg 0$, and therefore given any rational
$c > 0$ the multiplier ideal
$\MI{\tfrac{c}{k}| kL | }$ is defined for large $k$.
One checks easily that 
\begin{equation} \MI{\tfrac{c}{k} \cdot |kL| } \subseteq
\MI{\tfrac{c}{pk}\cdot |pkL|}  
\tag{*} \end{equation} for every integer $p > 0$. We
assert that then the family of ideals
$\big \{\MI{\tfrac{c }{k} \cdot |kL| } \big \}$ ($k \gg
0$) has a \textit{unique} maximal element. In fact, the
existence of at least one maximal member follows from
the ascending chain condition on ideals. On
the other hand, if
$\MI{\tfrac{c}{k } \cdot | k  L|}$ and
$\MI{\tfrac{c}{\ell} \cdot | \ell L|}$ are each maximal,
then thanks to (*)  they  must both coincide with
$\MI{\tfrac{c}{k\ell} \cdot | (k \ell) L|}$.
\begin{definition} The \textit{asymptotic multiplier
ideal sheaf} associated to $c$ and $|L|$, 
\[ \MI{ \, c \cdot \Vert L \Vert \, } \ = \ \MI{ \, X \
, \ c
\cdot \Vert L \Vert \, } ,
\] is defined to be the unique maximal member of the
family of ideals $\big \{ \, \MI{\tfrac{c }{k} \cdot
|kL| } \, \big \}$ ($k$~large). 
\end{definition}
\noi One can show that there exists a positive integer
$k_0$ such that $\MI{ c \cdot \Vert L \Vert} = \MI{
\tfrac{c}{k} \cdot |kL| }$ for every $k \ge k_0$. It
follows easily from the definition that $\MI{ m \cdot
\Vert L \Vert} = \MI{  \Vert mL \Vert}$ for every
positive integer $m > 0$.\footnote{In fact, fix $m > 0$.
Then for $p \gg 0$:
\[ \MI{ \Vert mL \Vert} = \MI{ \tfrac{1}{p} \cdot  
\vert mpL
\vert} = \MI{ \tfrac{m}{mp} \cdot \vert mpL \vert} =
\MI{ m
\cdot \Vert L \Vert}. \]} 

The basic facts about these asymptotic multiplier ideals
are summarized in the following
\begin{theorem} \label{properties.asymptotic.MI}  Let 
$X$ be a non-singular complex projective variety of
dimension $n$, and let $L$ be a big line bundle on
$X$.  
\begin{enumerate}
\item[(i)]  The natural inclusion
\[ \HH{0}{X}{\OO_X(L) \otimes \MI{\Vert L \Vert}} \lra
\HH{0}{X}{\OO_X(L)} \] is an isomorphism, i.e. 
$\MI{\Vert L \Vert }$ contains the base ideal ${\frak
b} \big( \vert L \vert \big )
\subset \OO_X$ of the linear series $\vert L \vert$. 
\item[(ii)] For any nef and big divisor $P$ one  has
the vanishing
\[ \HH{i}{X}{\OO_X(K_X + L + P) \otimes \MI{\Vert L
\Vert}} = 0 \for i > 0. \]
\item[(iii)]  If $B$ is very ample, then $\OO_X(K_X + 
L + (n+1)B) \otimes \MI{\Vert L \Vert}$ is generated by
its global sections. 
\end{enumerate}
\end{theorem}
\noi Of course the analogous statements hold  with $L$
replaced by $mL$.

\begin{proof}    The first statement follows easily
from the definition. For (ii) and (iii), note that
$\MI{\Vert L \Vert} = \MI{D}$ for a suitable
$\QQ$-divisor $D$ numerically equivalent to
$L$.  This being said,  (ii) is a consequence  of the
Nadel Vanishing theorem whereas (iii) follows from
Corollary \ref{MI.Gl.Gen}. \end{proof}

\begin{remark} The definition of the asymptotic
multiplier ideal
$\MI{\Vert L \Vert}$ requires only that $\kappa(X, L)
\ge 0$,  $\kappa(X, L)$ being the Kodaira-Iitaka
dimension of $L$, and Theorem
\ref{properties.asymptotic.MI} remains true in this
setting. When $L$ is big --- as we assumed for
simplicity --- the proof of Nadel Vanishing shows that
it suffices in statement (ii) that $P$ be nef, and
hence in (iii) one can replace the factor $(n+1)$ by
$n$. However we do not need these improvements here.
\qed
\end{remark}

Finally we discuss the relation between these
asymptotic multiplier ideals and their analytic
counterparts. In the analytic setting, there is a
concept of singular hermitian  metric
$h_{\min}$ with minimal singularities (see e.g.\
\cite{Demailly.Bourbaki}), defined whenever the first 
Chern class $c_1(L)$ lies in the closure of the cone
of effective divisors (``pseudoeffective cone'');
it is therefore not even necessary that
$\kappa(X,L)\ge 0$ for $h_{\min}$ to be defined, but only
that $L$ is pseudoeffective. The metric $h_{\min}$ is 
defined by taking any smooth hermitian metric $h_\infty$ 
on $L$ and putting $h_{\min}=h_\infty e^{-\psi_{\max}}$ where
\[
\psi_{\max}(x)=\sup\big\{\psi(x)\,;\,\psi~\mbox{usc},~\psi\le 0,~
i(\partial\overline\partial \log h_\infty+\psi)
\ge 0\big\}.
\]
For arbitrary sections $\sigma_1,{\ldots},\sigma_N\in H^0(X,kL)$
we can take \mbox{$\psi(x)={1\over k}\log\sum_j\Vert\sigma_j(x)
\Vert^2_{h_\infty}-C$} as an admissible $\psi$ function. 
We infer from this that the associated multiplier ideal
sheaf $\MI{h_{\min}}$ satisfies the inclusion 
\begin{equation} \label{Alg.vs.Analy.MI} \MI{\Vert
L\Vert}\subseteq 
\MI{h_{\min}}\end{equation}
  when $\kappa(X,L)\ge 0$.
The inclusion is strict in general. In fact, let us take
$E$ to be a unitary flat vector bundle on a smooth
variety
$C$ such that no non trivial symmetric power of $E$ or
$E^\star$ has sections  (such vector bundles exist 
already when $C$ is a curve of genus $\ge 1$), and
set $U = \OO_C \oplus E$. We take as our example $X =
\PP(U)$ and
$L=\OO_{\PP(U)}(1)$. Then for every $m \ge
1$, $\OO_X(mL)$ has a unique  nontrivial section which
vanishes to order
$m$ along the ``divisor at infinity" $H \subset \PP(U) =
X$, and hence $\MI{\Vert L\Vert}=\OO_X(-H)$. However
$L$ has a smooth semipositive metric induced by the
flat metric of $E$, so that $\MI{h_{\min}}=\OO_X$. It
is somewhat strange (but very interesting) that the
analytic setting yields ``virtual sections'' that do
not have algebraic counterparts. 

Note that in the example just presented, the
line bundle $L$ has Iitaka dimenson zero. We
conjecture that if $L$ is big, then equality should
hold in (\ref{Alg.vs.Analy.MI}).  We will 
prove here a slightly weaker statement, by
means of an analytic analogue of Theorem
\ref{properties.asymptotic.MI}. If $\phi$ is a 
plurisubharmonic function, the ideal sheaves
$\MI{(1+\eps)\phi}$ increase as $\eps$ decreases to
$0$, hence there must be a maximal element which we
denote by~$\MIP{\phi}$. This ideal always satisfies
$\MIP{\phi}\subseteq\MI{\phi}$. When $\phi$ has
algebraic singularities, standard semicontinuity
arguments show that $\MIP{\phi}=\MI{\phi}$, but we do
not know if equality  always holds in the analytic case.

\begin{theorem} Let $X$ be a non-singular complex projective variety of
dimension $n$, and let $L$ be a pseudoeffective line bundle 
on $X$.\footnote{Recall that the pseudoeffctive cone is
the closure of the cone of effective divisors on $X$.} 
Fix a singular hermitian metric
$h$ on
$L$ with nonnegative curvature current.
\begin{enumerate} \label{analytic.nadel}
\item[(i)] For any big and nef divisor $P$, one has the vanishing
\[ \HH{i}{X}{\OO_X(K_X + L + P) \otimes \MIP{h}} = 0 \for i > 0. \]
\item[(ii)]  If $B$ is very ample, then the sheaves
$\OO_X(K_X + L + (n+1)B) \otimes \MI{h}$ and\break $\OO_X(K_X + 
L + (n+1)B) \otimes \MIP{h}$ are generated by their global sections. 
\end{enumerate}
\end{theorem}

\begin{proof} (i) is a slight variation of Nadel's vanishing theorem
in its analytic form. If $P$ is ample, the result is true with
$\MI{h}$ as well as with $\MIP{h}$ (the latter case being obtained by
replacing $h$ with $h^{1+\eps}\otimes h_\infty^{-\eps}$ where $h_\infty$
is an arbitrary smooth metric on $L$; the defect of positivity of 
$h_\infty$ can be compensated by the strict positivity of $P$). If $P$ 
is big and nef, 
we can write $P=A+E$ with an ample $\QQ$-divisor $A$ and an effective 
$\QQ$-divisor $E$, and $E$ can be taken arbitrarily small. We then get 
vanishing with $\MIP{h\otimes h_E}$ where $h_E$ is the singular metric
of curvature current $[E]$ on $E$. However, if $E$ is so small that
$\MI{h_E^N}=\OO_X$, $N\gg 1$, we do have $\MIP{h\otimes h_E}=\MIP{h}$, 
as follows from an elementary argument using H\"older's inequality.

\noindent
Statement  (ii) follows from (i), Nadel Vanishing and
Mumford's Lemma \ref{CM.Reg}. Alternatively, one can
argue via   a straightforward adaptation of the proof
given  in
\cite{Siu.97}, based on Skoda's
$L^2$ estimates for ideals of holomorphic  functions
\cite{Skoda.72}.
\end{proof}  

\begin{theorem} Let $X$ be a projective nonsingular algebraic variety, $L$ 
a big nef line bundle on $X$, and $h_{\min}$ its singular hermitian metric
with minimal singularity. Then
\[
\MIP{h_{\min}}\subseteq\MI{\Vert L\Vert}\subseteq\MI{h_{\min}}.
\]
\end{theorem}

\begin{proof} The strong version of the Ohsawa-Takegoshi $L^2$ 
extension theorem proved by Manivel \cite{Manivel.93} shows that 
for every singular hermitian line bundle $(L,h)$ with nonnegative
curvature and every smooth complete intersection subvariety $Y\subseteq X$
(actually, the hypothesis that $Y$ is a complete intersection 
could probably be removed), there exists a sufficiently ample
line bundle $B$ and a surjective restriction morphism
\[
H^0\big(X,\OO_X(L+B)\otimes\MI{h}\big)
\longrightarrow
H^0\big(Y,\OO_Y(L+B)\otimes\MI{h_{|Y}}\big)
\]
with the following additional property: for every section on $Y$, there
exists an extension satisfying an $L^2$ estimate with a constant depending 
only on $Y$ (hence, independent of $L$). We take $Y$ equal to a smooth 
zero dimensional
scheme obtained as a complete intersection of hyperplane sections of
a very ample linear system $|A|$, and observe that $B$ depends only 
on $A$ in that case (hence can be taken independent of the choice of 
the particular $0$-dimensional scheme). Fix an integer $k_0$ so large that
$E:=k_0L-B$ is effective. We apply the extension theorem to the line
bundle $L'=(k-k_0)L+E$ equipped with the hermitian metric 
$h_{\min}^{k-k_0}\otimes h_E$, curv$(h_E)=[E]$ (and a smooth metric 
$h_B$ of positive curvature on $B$). Then, for $k\ge k_0$ and a prescribed
point $x\in X$, we select a zero-dimensional subscheme $Y$ containing $x$
and in this way we get a global section
$\sigma_x$ of $H^0(X,kL)=H^0(X,L'+E+B)$ such that
\[
\int_X\Vert\sigma_x(z)\Vert^2_{h_{\min}^{k-k_0}\otimes h_E\otimes h_B}
\le C\qquad
\mbox{while}\quad 
\Vert\sigma_x(x)\Vert_{h_{\min}^{k-k_0}\otimes h_E\otimes h_B}=1.
\] 
>From this we infer that locally $h_{\min}=e^{-2\phi}$
with
$|\sigma_x(x)|^2e^{-2(k-k_0)\phi(x)+2\phi_E+O(1)}=1$, hence
\[
\phi(x)+{1\over k-k_0}\phi_E\le
{1\over k-k_0}\log|\sigma_x(x)|+C\le
{1\over k-k_0}\log\sum_j|g_j(x)|+C
\]
where $(g_j)$ is an orthonormal basis of sections of $H^0(X,kL)$.
This implies that $\MI{\Vert h\Vert}$ contains the ideal
$\MI{h_{\min}\otimes h_E^{1/(k-k_0)}}$. Again, H\"older's inequality
shows that this ideal contains $\MIP{h_{\min}}$ for $k$ large enough.
\end{proof}

\section{Subadditivity}

The present section is devoted to the subadditivity
theorem stated in the Introduction, and some variants.

Let $X_1$, $X_2$ be smooth complex quasi-projective
varieties, and let $D_1$ and $D_2$ be effective
$\QQ$-divisors on
$X_1$, $X_2$, respectively.  Fix a log resolution 
$\mu_i : X_i^{\pr} \lra X_i$ of the pair $(X_i,D_i)$,
$i=1,2$.  We consider the product diagram
\[
\xymatrix{ X_1^{\pr} \ar[d]_{\mu_1} & X_1^\pr \times
X_2^{\pr}
\ar[l]_{q_1}
\ar[r]^{q_2}
\ar[d]^{\mu_1 \times \mu_2}  & X_2^{\pr} \ar[d]^{\mu_2}
\\ X_1 & X_1 \times X_2 \ar[l]^{p_1} \ar[r]_{p_2} & X_2
}
\] where the horizontal maps are projections.

\begin{lemma} The product $\mu_1 \times \mu_2 :
X_1^{\pr} \times X_2^{\pr}
\lra X_1 \times X_2$ is a log resolution of the pair
\[ ( \, X_1 \times X_2 \, , \,  p_1^*D_1 + p_2^* D_2 \,
). \]
\end{lemma}
\begin{proof}  Since the exceptional set $\exc( \mu_1
\times \mu_2)$ is the divisor where the derivative
$d(\mu_1 \times \mu_2)$ drops rank, one sees that 
$\exc(\mu_1 \times \mu_2) = q_1^* \exc(\mu_1) + q_2^*
\exc(\mu_2)$. Similarly, 
\[ (\mu_1 \times \mu_2)^* (p_1^*D_1 + p_2^* D_2) \ = \
q_1^* \mu_1^* D_1 + q_2^* \mu_2^* D_2. \]  Therefore
\small
\[ \exc(\mu_1 \times \mu_2) + (\mu_1 \times \mu_2)^*
(p_1^*D_1 + p_2^* D_2) \ = \ q_1^*(\exc(\mu_1) +
\mu_1^* D_1) + q_2^*(\exc(\mu_2) + \mu_2^* D_2), \]
\normalsize and this has normal crossing support since
$\exc(\mu_1) +
 \mu_1^* D_1$ and $\exc(\mu_2) +  \mu_2^* D_2$ do. 
\end{proof}

\begin{proposition} \label{diag.product} One has
\[ \MI{ \, X_1 \times X_2 \, , \, p_1^* D_1 + p_2^*
D_2} \ = \ p_1^{-1} \MI{X_1, D_1} \cdot p_2^{-1}
\MI{X_2, D_2}. \]
\end{proposition}
\begin{proof} To lighten notation we will write $D_1
\boxplus D_2$ for the exterior direct sum $p_1^* D_1 +
p_2^* D_2$, so that the formula to be established is
\[ 
\MI{ \, X_1 \times X_2 \, , \,  D_1 \boxplus   D_2 \,} \
= \ p_1^{-1} \MI{X_1, D_1} \cdot p_2^{-1} \MI{X_2, D_2}.
\]  The plan is to compute the  multiplier ideal on the
left using the log resolution
$\mu_1 \times \mu_2$. Specifically:
\[
\MI{ \, X_1 \times X_2 \, , \,  D_1 \boxplus   D_2 \, }
\ =
\ (\mu_1 \times \mu_2 )_*   \OO_{X_1^\pr \times
X_2^\pr}\big( K_{X_1^\pr \times X_2^\pr / X_1  \times
X_2 } - [ (\mu_1
\times \mu_2)^* ( D_1 \boxplus D_2)]  \big ). \] Note
to begin with that
\[
\big [ (\mu_1 \times \mu_2)^* (D_1 \boxplus D_2) \big ]
=  \big [ q_1^* \mu_1^* D_1 \big ] \, + \, \big [ q_2^*
\mu_2^* D_2 \big ] 
\] thanks to the fact that $q_1^* \mu_1^* D_1$ and
$q_2^*
\mu_2^* D_2$ have no common components. Furthermore, as
$q_1$ and $q_2$ are smooth:
\[
 \ \big[ q_1^* \mu_1^* D_1 \big] \, = \, q_1^* \big[
  \mu_1^* D_1
\big]  \ \quad \text{and}\ \quad  
\big [ q_2^*  \mu_2^* D_2 \big]\, = \, q_2^*  \big[  
\mu_2^* D_2 \big].  \]
 Since 
$ K_{X_1^\pr \times X_2^\pr / X_1 \times X_2 } \ = \
q_1^* \big( K_{X_1^\pr / X_1} \big ) \, + \, q_2^* \big(
K_{X_2^\pr / X_2} \big )$, it then follows that 
\begin{multline*}
 \OO_{X_1^\pr \times X_2^\pr}\Big( K_{X_1^\pr\times
X_2^\pr / X_1  \times X_2 } - [ (\mu_1 \times
\mu_2)^* (p_1^* D_1 + p_2^* D_2)]  \Big ) \\ = \    
q_1^*
\OO_{X_1^\pr}\big( K_{X_1^\pr / X_1} - [ \mu_1^* D_1]
\big)
\otimes  q_2^* \OO_{X_2^\pr}\big( K_{X_2^\pr / X_2} - [
\mu_2^* D_2] \big). 
\end{multline*} Therefore
\begin{align*}
 \MI{ \, X_1 \times X_2 \, , \, D_1 \boxplus D_2 \, }
\kern-40pt&\kern40pt = \\ &= \ (\mu_1 \times \mu_2)_*
\Big( q_1^*
\OO_{X_1^\pr}\big( K_{X_1^\pr / X_1} - [ \mu_1^* D_1]
\big)
\otimes  q_2^* \OO_{X_2^\pr}\big( K_{X_2^\pr / X_2} - [
\mu_2^* D_2] \big) 
\Big ) \\ &= \  p_1^* \mu_{1\,*} \OO_{X_1^\pr}\big(
K_{X_1^\pr / X_1} - [
\mu_1^* D_1] \big) \
\otimes  \ p_2^* \mu_{2\,*} \OO_{X_2^\pr}\big(
K_{X_2^\pr / X_2} - [ \mu_2^* D_2] \big)  \\ &= \ p_1^*
\MI{X, D_1}\, \otimes
\,  p_2^* \MI{X, D_2}
\end{align*} thanks to the K\"unneth formula. But
\[ p_1^* \MI{X_1, D_1} = p_1^{-1} \MI{X_1, D_1} \ \
\text{and} \ \ p_2^* \MI{X_2, D_2} = p_2^{-1} \MI{X_2,
D_2}
\] since $p_1$ and $p_2$ are flat. Finally, 
\[ p_1^{-1} \MI{X_1, D_1}\, \otimes
\,  p_2^{-1} \MI{X_2, D_2}  \ = \ p_1^{-1} \MI{X_1,
D_1}\,
\cdot \, p_2^{-1} \MI{X_2, D_2} \] by virtue of the
fact that $p_1^{-1}\MI{X_1, D_1}$ is flat for  $p_2$
(cf.\ \cite{Matsumura}). This completes the proof of the
Proposition. 
\end{proof}

The subadditivity property of multiplier ideals now
follows immediately:
\begin{theorem} \label{Subadditivity.Thm} Let $X$ be a
smooth complex quasi-projective variety, and let $D_1$
and $D_2$ be effective $\QQ$-divisors on
$X$. Then
\[ \MI{ X, D_1 + D_2} \, \subseteq \, \MI{X,D_1} \cdot
\MI{X, D_2}. \]
\end{theorem}
\begin{proof} We apply Corollary
\ref{restr.to.submnfld} to the diagonal $ \Delta = X 
\subset X\times X$. Keeping the notation of the previous
proof (with $X_1=X_2=X$, $\mu_1=\mu_2=\mu$), one has
\begin{align*}
\MI{X, D_1 + D_2} \ &= \MI{ \, \Delta \, , \,  \big(
p_1^* D_1 +
 p_2^* D_2 \big)_{|\Delta} } \\ &\subseteq \MI{ \,  X
\times X \, , \, p_1^* D_1 + p_2^* D_2 \, \, } \cdot
\OO_{\Delta}
\end{align*} But it follows from Proposition
\ref{diag.product} that
\[ \MI{ \,  X \times X \, , \, p_1^* D_1 + p_2^* D_2 \,
\, } \cdot \OO_{\Delta} \ = \ \MI{X, D_1} \cdot
\MI{X,D_2}, \] as required. 
\end{proof}

\begin{variant} \label{subadditivity.asmypt.ideals} Let
$L$ be a big line bundle on a non-singular complex
projective variety $X$. Then for all $m \ge 0$:
\[ \MI{X, \Vert mL \Vert} \subseteq \MI{X, \Vert L
\Vert}^m. \]
\end{variant}
\begin{proof} Given $m$, fix $p \gg 0$ plus a general
divisor $D \in
\vert mpL \vert$. Then
\[ \MI{\Vert L \Vert} \ = \ \MI{ \tfrac{1}{pm}D} \ \
\text{and} 
\ \ \MI{\Vert mL \Vert} \ = \ \MI{\tfrac{1}{p} D}, \]
so the assertion follows from the Theorem. 
\end{proof}

\begin{variant} Let $\fra , {\frak b} \subseteq \OO_X$
be ideals, and fix rational numbers $c , d > 0$. Then
\[ \JJ \big( (c   \cdot  \fra ) \cdot ( d \cdot \frb )\big) 
\ \subseteq \ 
\MI{ c \cdot \fra} \cdot \MI{d \cdot {\frak b}}. \]
\end{variant}
\begin{proof} This does not follow directly from the
statement of Theorem \ref{Subadditivity.Thm} because
the divisor of a general element of $\fra \cdot {\frak
b}$ is not the sum of divisors of elements in $\fra$
and ${\frak  b}$. However the proof Proposition
\ref{diag.product} goes through  to show that
\[
\JJ \Big(\,  X \times X \, , \, \big(c   \cdot  p_1^{-1}
\fra \big ) 
\cdot \big ( d \cdot p_2^{-1} {\frak b} \big) \, \Big) \
=
\ p_1^{-1}
\MI{X, c \cdot \fra} \cdot p_2^{-1} \MI{X , d \cdot
{\frak b}}, \] and then as above one restricts to the
diagonal. \end{proof}

The corresponding properties of   analytic multiplier
ideals are proven in the analogous manner. The
result is the following: 
\begin{theorem} [Analogous analytic statements]
\begin{enumerate}
\item[]
\item[(i)] Let $X_1$, $X_2$  be complex manifolds and
let $\phi_i$ be a plurisubharmonic function on $X_i$.
Then
\[
\MI{\phi_1\circ p_1 + \phi_2\circ p_2} = p_1^{-1}
\MI{\phi_1} \cdot  p_2^{-1}\MI{\phi_2}.
\]
\item[(ii)] Let $X$ be a complex manifold and let
$\phi$, $\psi$ be plurisubharmonic functions on $X$.
Then
\[
\MI{\phi + \psi} \subseteq \MI{\phi} \cdot \MI{\psi}
\]
\end{enumerate}
\end{theorem}
\begin{proof} Only (i) requires a proof, since (ii)
follows again from (i)  by the restriction principle
and the diagonal trick. Let us fix two relatively
compact Stein open subsets $U_1\subset X_1$, 
$U_2\subset X_2$. Then $\cH^2(U_1\times U_2,
\phi_1\circ p_1 + \phi_2\circ p_2, p_1^\star
dV_1\otimes p_2^\star dV_2)$ is the Hilbert tensor
product of $\cH^2(U_1,\phi_1,dV_1)$ and
$\cH^2(U_2,\phi_2,dV_2)$, and admits $(f'_k\boxtimes
f''_l)$ as a Hilbert  basis, where $(f'_k)$ and
$(f''_l)$ are respective Hilbert bases. Since
$\MI{\phi_1\circ p_1 + \phi_2\circ p_2}_{|U_1\times
U_2}$ is generated  as an $\OO_{U_1\times U_2}$ module
by the $(f'_k\boxtimes f''_l)$, we  conclude that (i)
holds true.
\end{proof}

\section{Fujita's Theorem}

Now let $X$ be a smooth irreducible complex projective
variety of dimension $n$, and $L$ a line bundle on
$X$. We recall the
\begin{definition} The \textbf{volume}\footnote{This
was called the ``degree" of the graded linear series
$\oplus \HH{0}{X}{\OO_X(kL)}$ in
\cite{ELN}, but the present terminology is more natural
and seems to be becoming standard.} of $L$ is the real
number
\[ v(L) \ = \ v(X, L) \ = \ \limsup_{k \to \infty} \, 
\frac{n!}{k^n} 
 \, \hh{0}{X}{\OO(kL)}. \qed \]
\end{definition}
\noi  Thus $L$ is big iff
$v(L) > 0$. If $L$ is ample, or merely nef and big,
then asymptotic Riemann-Roch shows that
\[
\hh{0}{X}{\OO_X(kL)} = \frac{k^n}{n!} \big( L^n \big )
+ o(k^n),
\] so that in this case $v(L) = \big ( L^n \big)$ is the
top self-intersection number of $L$. If $D$ is a
$\QQ$-divisor on $X$, then the volume $v(D)$ is defined
analogously, the limit being taken over $k$ such that
$kD$ is an integral divisor. 

Fujita's Theorem asserts that ``most of" the volume of
$L$ can be accounted for by the volume an ample
$\QQ$-divisor on a modification.
\begin{theorem} [Fujita  \cite{Fujita}] 
\label{Fujita.Thm}  Let
$L$ be a big line bundle on $X$, and fix 
$\eps > 0$. Then there exists a birational modification
\[
\mu : X^\pr      \lra X \] (depending on $\epsilon$)
and a decomposition $
\mu^* L \lin E  + A $ (also depending on $\epsilon$), 
with
$E $  an effective $\QQ$-divisor and $A$ an  ample
$\QQ$-divisor on $S^{\pr}$, such that \[ v(X^\pr, A) =
\big( A^n \big) \ge v(X, L) -
\epsilon.\]  \end{theorem}
\noi Conversely, given a decomposition $\mu^* L   \lin E
+ A$ as in the Theorem, one evidently has $v(X^\pr , A)
=
\big( A^n ) \le v(X, L)$. So the essential content of
Fujita's theorem is that the volume of any big line
bundle can be approximated arbitrarily closely by the
volume of an ample $\QQ$-divisor (on a modification).
This statement initially arose in connection with
alegbro-geometric analogues of the work
\cite{Demailly.93} of the first author (cf. \cite{LLS},
\S 7; \cite{ELN}). A geometric reinterpretation appears
in Proposition  \ref{moving.self.int}. 

\begin{remark} Suppose that $L$ admits a \textit{Zariski
decomposition}, i.e. assume that there exists a
birational modification 
$\mu : X^{\pr} \lra X$, plus a decomposition $\mu^* L =
P + N$, where $P$ and $N$ are $\QQ$-divisors, with $P$
nef, having the property that \[ \HH{0}{X}{\OO_X(kL)} =
\HH{0}{X^\pr}{\OO_{X^{\pr}}( [kP] )}\]  for all $k \ge
0$. Then $v(X, L) = v(X^\pr,  P) = \big( P^n \big)$,
i.e. the volume of
$L$ is computed by the volume of a nef divisor on a
modification.  While is it known that such
decompositions do not exist in general
\cite{Cutkosky}, Fujita's Theorem shows that an
approximate asymptotic statement does hold. \qed
\end{remark}

Fujita's proof is quite short, but rather tricky: it is
an argument by contradiction revolving around the Hodge
index theorem. The purpose of this section is to use
the subadditivity property of multiplier ideals to give
a new proof which seems perhaps a bit more transparent.
(One can to a certain extent see the present argument as
extending to all dimensions the proof for surfaces due
to Fernandez del Busto appearing in
\cite{LLS}, \S 7.)

We begin with two lemmas. The first, due to Kodaira, is
a standard consequence of asymptotic Riemann-Roch (cf.
\cite{Kollar.Book}, (VI.2.16)).
\begin{lemma} [Kodaira's Lemma]
\label{Kodaira.Lemma} Given a big line bundle $L$, and
any ample bundle $A$ on~$X$, there is a positive
integer $m_0 > 0$ such that $m_0 L = A + E $ for some
effective divisor~$E$. \qed
\end{lemma}

The second (somewhat technical) Lemma shows that one can
perturb
$L$ slightly without greatly affecting its volume:
\begin{lemma} \label{Volume.Continuous} Let $G$ be an
arbitrary line bundle. For every
$\eps >0$, there exists a positive integer $m$ and a
sequence
$\ell_\nu\uparrow+\infty$ such that
\[
\hh{0}{X}{\ell_\nu(mL-G)} \ge {\ell_\nu^n m^n\over n!}
\big(v(L)-\eps \big).
\] In other words, \[ v(mL-G) \ \ge\  m^n
\big( v(L)-\epsilon \big) \] for
$m$ sufficiently large.
\end{lemma}

\begin{proof} Clearly, $v(mL-G)\ge v(mL-(G+E))$ for
every effective divisor
$E$. We can take $E$ so large that $G+E$ is very ample,
and we are thus  reduced to the case where $G$ itself is
very ample by replacing $G$ with $G+E$. By definition of
$v(L)$, there exists a sequence $k_\nu\uparrow+\infty$
such that
\[ \hh{0}{X}{\OO_X(k_\nu L)} \ \ge \ {k_\nu^n\over
n!}\Big(v(L)-\frac{\eps}{2} \Big).
\] We now fix an integer $m\gg 1$ (to be chosen
precisely  later), and put
$\ell_\nu=
\big[{k_\nu\over m}\big]$, so that $k_\nu=\ell_\nu
m+r_\nu$,
$0\le r_\nu<m$. Then
\[
\ell_\nu(mL-G)=k_\nu L-(r_\nu L+\ell_\nu G).
\] Fix next a constant $a\in {\mathbf N} $ such that
$aG-rL$ is an effective divisor for each $0 \le r < m$.
Then
$maG - r_\nu L$ is effective, and  hence
\[ \hh{0}{X}{\OO_X(\ell_\nu(mL-G))} \ \ge \
\hh{0}{X}{\OO_X(k_\nu L-(\ell_\nu+am)G)}. \]
  We select a smooth divisor $D$ in the very ample
linear system
$|G|$. By looking at global sections associated with
the exact sequences of sheaves
\small
\[ 0 \lra \OO_X(-(j+1)D)\otimes \OO_X(k_\nu L) \lra
\OO_X(-jD)\otimes\OO_X(k_\nu L)
\lra \OO_D(k_\nu L-jD)\lra 0,
\]
\normalsize
$0\le j<s$, we infer inductively that
\begin{align*} 
\hh{0}{X}{\OO_X(k_\nu L-sD)} &\ge
\hh{0}{X}{\OO_X(k_\nu L)}-\sum_{0\le j<
s}\hh{0}{D}{\OO_D(k_\nu L-jD)} 
\\  &\ge\hh{0}{X}{\OO_X(k_\nu L)} -
s\,\hh{0}{D}{\OO_D(k_\nu L)} \\
 &\ge \frac{k_\nu^n}{n!}\Big(v(L)-{\frac{\eps}{2}}\Big)
-s\,Ck_\nu^{n-1}
\end{align*}
 where $C$ depends only on $L$ and~$G$. Hence, by
putting $s=\ell_\nu+am$, we get
\begin{align*}
 \hh{0}{X}{\OO_X(\ell_\nu(mL-G))}   &\ge
\frac{k_\nu^n}{n!}\Big( v(L)-\frac{\eps}{2}\Big) 
-C(\ell_\nu+am)k_\nu^{n-1}\\
 &\ge \frac{\ell_\nu^nm^n}{
n!}\Big(v(L)-\frac{\eps}{2}\Big)
-C(\ell_\nu+am)(\ell_\nu+1)^{n-1}m^{n-1}
\end{align*}  and the desired conclusion follows by
taking
$\ell_\nu\gg m\gg 1$. 
\end{proof}

Now we turn to the 

\begin{proof} [Proof of Fujita's Theorem] Note to begin
with that it is enough to produce a big and nef divisor
$A$ satisfying the conclusion of the Theorem. For by
Kodaira's Lemma one can write $A \lin E + A^\pr$
where
$E$ is an effective
$\QQ$-divisor, and $A^\pr$ is an ample $\QQ$-divisor.
Then 
\[ E + A \ \lin \ E + \delta E  + (1 - \delta) A +
\delta A^{\pr}, \] where $A^{\pr \pr} =_{\text{def}} (1
- \delta)A +
\delta A^\pr$ is ample and the top self intersection
number $\big( (A^{\pr \pr})^n \big)$ approaches $\big (
A^n \big)$ as closely as we want.  

Fix now a very ample bundle $B$ on $X$,  set $G = K_X +
(n+1)B$, and for $m \ge 0$ put
\[ M_m = mL - G. \]  We can suppose that $G$ is very
ample, and we choose a  divisor
$D \in |G|$.  Then multiplication by $\ell D$ determines
for every $\ell \ge 0$ an inclusion $\OO_X(\ell M_m)
\hookrightarrow \OO_X(\ell mL)$  of sheaves, and
therefore an injection
\[\HH{0}{X}{\OO_X(
\ell M_m)}
\subseteq \HH{0}{X}{\OO_X(\ell mL)}.\]   Given
$\epsilon > 0$, we use Lemma 
\ref{Volume.Continuous} to fix  $m \gg 0$ such that
\begin{equation} v(M_m)
\ge m^n \big (v(L) - \epsilon \big). 
\label{large.m} \end{equation}  We further assume  that
$m$ is sufficiently large so that
$M_m$ is big. 

Having fixed $m
\gg 0$ satisfying (\ref{large.m}), we will produce an
ideal sheaf
$\JJ = \JJ_m \subset \OO_X$ (depending on
$m$) such that 
\begin{gather}  
\OO_X(mL) \otimes \JJ   \ \text{is globally generated;}
\label{property.i} \\
\HH{0}{X}{\OO_X(\ell M_m)} \subseteq
\HH{0}{X}{\OO_X(\ell m L) \otimes \JJ^\ell} \fall
\ell \ge 1. \label{property.ii} 
\end{gather} Granting for the time being the existence
of $\JJ$, we complete the proof. Let $\mu : X^\pr \lra
X$ be a log resolution of $\JJ$, so that $\mu^{-1} \JJ =
\OO_{X^{\pr}}(-E_m)$ for some effective divisor $E_m$ on
$X^\pr$. It follows from (\ref{property.i}) that 
\[ A_m \ =_{\text{def}} \ \mu^* \big(mL \big)  - E_m \]
is globally generated, and hence nef. Using
(\ref{property.ii}) we find:
\begin{align*}
\HH{0}{X}{\OO_X(\ell M_m)} &\subseteq
\HH{0}{X}{\OO_X(\ell m L) \otimes \JJ^\ell} \\ 
&\subseteq \HH{0}{X^\pr}{\OO_{X^{\pr}}\big(\mu^{*}(\ell
m L) - \ell E_m \big)} \\ &=
\HH{0}{X^\pr}{\OO_{X^\pr}(\ell A_m)} 
\end{align*} (which shows in particular that $A_m$ is
big). This implies that 
\begin{align*}
\big( (A_m)^n \big) \ &= \ v(X^\pr, A_m) \\
        & \ge v(X, M_m) \\ & \ge m^n \big( v( L) -
\epsilon \big),
\end{align*} so the Theorem follows upon setting $A =
\tfrac{1}{m} A_m$ and $E = \tfrac{1}{m} E_m$. 

Turning to the construction of $\JJ$, set
\[ \JJ = \MI{X, \Vert M_m \Vert }. \] Since $mL = M_m +
\big(K_X + (n+1)B \big)$, (\ref{property.i}) follows
from Theorem
\ref{properties.asymptotic.MI}(iii) applied to $M_m$. 
As for (\ref{property.ii})  we  first apply Theorem
\ref{properties.asymptotic.MI}(i) to $ \ell M_m $,
together with the subadditivity property  in the form of
Variant 
\ref{subadditivity.asmypt.ideals}, to conclude:
\begin{equation} \label{first.eqn}
\begin{aligned}
 \HH{0}{X}{\OO_X(\ell M_m)} \ &=  \
\HH{0}{X}{\OO_X( \ell M_m) \otimes \MI{    \Vert
\ell M_m \Vert} }\\ &\subseteq \ \HH{0}{X}{\OO_X(\ell
M_m) \otimes
\MI{\Vert M_m \Vert}^{\ell}}.
\end{aligned}
\end{equation} Now  the sheaf homomorphism
\[ \OO_X(\ell M_m ) \otimes \MI{\Vert M_m
\Vert}^\ell \overset{\cdot \ell D}{\lra} \OO_X(\ell m L)
\otimes \MI{\Vert M_m \Vert}^\ell \] evidently remains
injective for all $\ell$, and consequently
\begin{equation} \label{second.eqn}
\HH{0}{X}{\OO_X(\ell M_m) \otimes
\MI{\Vert M_m \Vert}^{\ell}}
\subseteq \ \HH{0}{X}{\OO_X(\ell m L) \otimes \
\MI{\Vert M_m \Vert}^{\ell}}.
\end{equation} The required inclusion
(\ref{property.ii}) follows by combining
(\ref{first.eqn}) and (\ref{second.eqn}). This
completes the proof of Fujita's Theorem.
\end{proof}

We conclude by using Fujita's theorem to establish a 
geometric interpretation of the volume $v(L)$. Suppose
as above that $X$ is a smooth projective variety of
dimension
$n$, and that $L$ is a big line bundle on $X$. Given a
large integer
$k \gg 0$, denote by $B_k \subseteq X$ the base-locus
of the linear series $| kL |$. The \textbf{moving
self-intersection number} $\big(kL\big)^{[n]}$ of $|kL|$
is defined by choosing $n$ general divisors $D_1, \dots
, D_n \in |kL|$ and putting
\[ \big( kL \big)^{[n]} \ = \ \# \Big( D_1 \cap \ldots
\cap D_n
\cap ( X - B_k)  \Big). \] In other words, we simply
count the number of intersection points away from the
base locus of $n$ general divisors in the linear series
$|kL|$. This notion arises for example in Matsusaka's
proof of his ``big theorem"
 (cf.\cite{LM}).

We show that the volume $v(L)$ of $L$ measures the rate
of growth with respect to $k$ of these moving
self-intersection numbers. The following result is
implicit in
\cite{Tsuji}, and was undoubtably known also to Fujita.
\begin{proposition} \label{moving.self.int} Assume as
above that $L$ is a big line bundle on a smooth
projective variety $X$. Then one has
\[ v(L) \ = \ \limsup_{k \to \infty} \frac{ \big( kL
\big)^{[n]} }{k^n}. \]
\end{proposition}
\begin{proof} We start by interpreting $\big(
kL\big)^{[n]}$ geometrically. Let $ \mu_k : X_k \lra X$
be a log resolution of $|kL|$, with
$\mu_k^* |kL| = |V_k| +  F_k $, where
 \[ P_k \  =_{\text{def}}  \ \mu_k^* (\, kL \, ) \, - \,
F_k\]  is free, and
$\HH{0}{X}{\OO_X(kL)} = V_k =
\HH{0}{X_k}{\OO_{X_k}(P_k)}$, so that $B_k =
\mu_k(F_k)$. Then evidently
$(kL)^{[n]}$ counts the number of intersection points
of $n$ general divisors in $P_k$, and consequently
\[ \big( kL \big)^{[n]} \  = \ 
\big( (P_k)^n \big).\]  We have 
$\big( (P_k)^n \big ) = v(X_k, P_k)$ for $k \gg 0$
since then $P_k$ is big (and nef), and $v(X, kL) \ge
v(X_k, P_k)$ since $P_k$ embeds in
$\mu_k^*(k L)$.  Hence
\[ v(X, kL) \ \ge  \ \big( kL
\big)^{[n]} \ \for k \gg 0. \] On the other hand, an
argument in the spirit of Lemma \ref{Volume.Continuous}
shows that $v(X, kL) = k^n \cdot v(X, L)$ (\cite{ELN},
Lemma 3.4), and so we conclude that 
\begin{equation} v(L) \ge \frac{ \big( kL
\big)^{[n]} }{k^n}. \tag{*}
\end{equation} for every $k \gg 0$. 

For the reverse inequality we use Fujita's theorem. Fix
$\epsilon > 0$, and consider the  decomposition $\mu^*
L = A + E$ on $\mu : X^\pr \lra X$ constructed in
(\ref{Fujita.Thm}). Let $k$ be any positive integer
such that
$kA$ is integral and globally generated. By taking a
common resolution we can assume that $X_k$ dominates
$X^\pr$, and hence we can write 
\[ \mu_k^* \, kL \ \lin \  A_k + E_k \] with $A_k$
globally generated and \[ \big( (A_k)^n \big)  \ \ge
\ k^n
\cdot \big( v(X,L) - \epsilon \big).\]  But then
$\HH{0}{X_k}{A_k}$ gives rise to  a free linear
subseries of
$\HH{0}{X_k}{P_k}$, and consequently \[
\big( (A_k)^n
\big) \ \le \ \big( (P_k)^n \big) \  = \ \big( kL
\big)^{[n]}.
\]
 Therefore
\begin{equation}
\frac{ \big( kL \big)^{[n]}}{k^n} \ \ge \ v(X,L) -
\epsilon. \tag{**}
\end{equation} But (**) holds for any sufficiently
large and divisible
$k$, and in view of (*) the Proposition follows. 
\end{proof}

\end{document}